\documentclass[12pt]{amsart}
\usepackage{amssymb}

\usepackage{mathptmx}

\usepackage{pstcol}
\psset{unit=1.0cm,linewidth=0.4pt}
\def\vertex{\pscircle[fillstyle=solid,fillcolor=black]{0.05}}
\definecolor{light}{gray}{0.9}
\definecolor{medium}{gray}{0.8}

\newcommand{\mm}{\mathfrak m}

\newcommand{\pp}{\mathfrak p}
\newcommand{\qq}{\mathfrak q}

\newcommand{\Z}{\mathbb{Z}}
\newcommand{\R}{\mathbb{R}}
\newcommand{\N}{\mathbb{N}}

\newcommand{\Fc}{\mathcal{F}}

\DeclareMathOperator{\chara}{char} \DeclareMathOperator{\Hom}{Hom}
 \DeclareMathOperator{\Ker}{Ker}
 \DeclareMathOperator{\rank}{rank}

\DeclareMathOperator{\tensor}{\otimes}
\DeclareMathOperator{\iso}{\cong}

\DeclareMathOperator{\pnt}{\raise 0.5mm \hbox{\large\bf.}}

\DeclareMathOperator{\gp}{G}
\DeclareMathOperator{\cn}{C}
\DeclareMathOperator{\conv}{conv}

\DeclareMathOperator{\depth}{depth}
\DeclareMathOperator{\Spec}{Spec}

\DeclareMathOperator{\Ext}{Ext} 

\DeclareMathOperator{\relint}{int} \DeclareMathOperator{\nat}{nat}

\def\+#1{\relax\ifmmode\if\noexpand #1\relax \mathop{\kern
    0pt^+{#1}}\nolimits\else \kern 0pt^+\!#1 \fi\else$^*$#1\fi}

\newtheorem{thm}{\bf Theorem}[section]
\newtheorem{lem}[thm]{\bf Lemma}
\newtheorem{cor}[thm]{\bf Corollary}
\newtheorem{prop}[thm]{\bf Proposition}

\theoremstyle{definition}

\newtheorem{rem}[thm]{\bf Remark}
\newtheorem{ex}[thm]{\bf Example}

\theoremstyle{plain}
\newtheorem*{thm*}{Theorem}

\title{On seminormal monoid rings}

\author{Winfried Bruns}
\address{FB Mathematik/Informatik, Universit\"at Osnabr\"uck, 49069 Osnabr\"uck, Germany}
\email{winfried@mathematik.uni-osnabrueck.de}

\author{Ping Li}
\address{
Dept. of Math. \& Stat., Queen's University, Kingston, On, K7L 3N6, Canada}
\email{pingli@mast.queensu.ca}

\author{Tim R\"omer}
\address{FB Mathematik/Informatik, Universit\"at Osnabr\"uck, 49069 Osnabr\"uck, Germany}
\email{troemer@mathematik.uni-osnabrueck.de}

\begin{document}

\begin{abstract}
Given a seminormal affine monoid $M$
we consider several monoid properties of $M$ and their connections to
ring properties of the associated affine monoid ring $K[M]$ over a field $K$.
We characterize when $K[M]$ satisfies Serre's condition ($S_2$)
and analyze the local cohomology of $K[M]$.
As an application we present
criteria which imply that $K[M]$ is
Cohen--Macaulay and we give lower bounds for the depth of $K[M]$.
Finally, the seminormality of an arbitrary affine monoid $M$
is studied with characteristic $p$ methods.
\end{abstract}

\maketitle \tableofcontents
%------------------------------------------------------------------------
%
%
%
%------------------------------------------------------------------------
\section{Introduction}
Let $M$ be an affine monoid, i.e. $M$ is a finitely generated
commutative monoid which can be embedded into $\Z^m$ for some $m \in
\N$. Let $K$ be a field and $K[M]$ be the affine monoid ring
associated to $M$. Sometimes $M$ is also called an affine semigroup
and $K[M]$ a semigroup ring. The study of affine monoids and affine
monoid rings has applications in many areas of mathematics.
It establishes
the combinatorial background for the theory of toric
varieties, which is the strongest connection to algebraic
geometry. In the last decades many authors have studied the
relationship between ring properties of $K[M]$ and monoid properties
of $M$. See Bruns and Herzog \cite{BRHE98} for a detailed discussion
and Bruns, Gubeladze and Trung \cite{BRGUTR02} for a survey about
open problems.

A remarkable result of Hochster \cite{HO72} states that if $M$ is a
normal, then $K[M]$ is Cohen--Macaulay. The converse is not true. It
is a natural question to characterize the Cohen--Macaulay property
of $K[M]$ for arbitrary affine monoids $M$ in terms of combinatorial
and topological information related to $M$. Goto, Suzuki and
Watanabe \cite{GSW76} could answer this question for simplicial
affine monoids. Later Trung and Hoa \cite{TRHO86} generalized their
result to arbitrary affine monoids. But the characterization is
technical and not easy to check. Thus it is interesting to consider
classes of monoids which are not necessarily simplicial, but
nevertheless admit simple criteria for the Cohen-Macaulay property.

One of the main topics in the thesis \cite{LI04} of the second
author were seminormal affine monoids and their monoid rings. Recall
that an affine monoid $M$ is called {\em seminormal} if
$z\in\gp(M)$, $z^2 \in M$ and $z^3 \in M$ imply that $z \in M$. Here
$\gp(M)$ denotes the group generated by $M$. Hochster and Roberts
\cite[Proposition 5.32]{HORO76} noted that $M$ is seminormal if and
only if $K[M]$ is a seminormal ring. See Traverso \cite{TR70} or
Swan \cite{SW80} for more details on the latter subject. In general
there exist Cohen--Macaulay affine monoid rings which are not
seminormal, and there exist seminormal affine monoid rings which are
not Cohen--Macaulay. One of the main goals of this paper is to
understand the problem in which cases $K[M]$ is Cohen--Macaulay for
a seminormal affine monoid $M$. Another question is to characterize
the seminormality property of affine monoids.

Let us go into more detail. Let $R$ be a Noetherian ring and let $N$
be a finitely generated $R$-module. The module $N$ satisfies {\em
Serre's condition} ($S_k$) if
$$
\depth N_\pp \geq \min\{k,\dim N_\pp\}
$$
for all $\pp \in \Spec R$. For trivial reasons Cohen--Macaulay rings
satisfy Serre's condition ($S_k$) for all $k \geq 1$. The main
result of \cite{GSW76} and a result of Sch\"afer and Schenzel
\cite{SS90} show that for a simplicial affine monoid $M$ the ring
$K[M]$ is Cohen--Macaulay if and only  if $K[M]$ satisfies ($S_2$).
After some prerequisites we study in Section \ref{s2section} the
question to characterize the ($S_2$) property for $K[M]$ if $M$ is a
seminormal monoid. In the following let $\cn(M)$ be the cone
generated by $M \subseteq \Z^m$. The main result in this section
already appeared in the thesis of the second author \cite{LI04} and
states:

\begin{thm*}
Let $M \subseteq \Z^m$ be a seminormal monoid
and let $F_1,\dots,F_t$ be the facets of
$\cn(M)$. Then the following statements are equivalent:
\begin{enumerate}
\item
$K[M]$ satisfies $(S_2)$;
\item
For all proper faces $F$ of $\cn(M)$ one has
$$
M \cap \relint F=\bigcap_{F \subseteq F_j} \gp(M\cap F_j)\cap
\relint F;
$$
\item
$\gp(M\cap F)=\bigcap_{F \subseteq F_j} \gp(M\cap F_j)$.
\end{enumerate}
\end{thm*}

Here $\relint F $ denotes the relative interior of $F$ with respect to the
subspace topology on the affine hull of $F$. Let us assume for a
moment that $M$ is positive, i.e.\ $0$ is the only invertible
element in $M$. In order to decide whether $K[M]$ is a
Cohen--Macaulay ring, one must understand the local cohomology
modules $H^i_\mm(K[M])$ where $\mm$ denotes the maximal ideal of
$K[M]$ generated by all monomials $X^a$ for $a \in M \setminus
\{0\}$, because the vanishing and non-vanishing of these modules
control the Cohen--Macaulayness of $K[M]$. Already Hochster and
Roberts \cite{HORO76} noticed that certain components of
$H^i_\mm(K[M])$ vanish for a seminormal monoid. Our result
\ref{semihelper2} in Section \ref{localcohom} generalizes their
observation, and we can prove the following

\begin{thm*}
Let $M\subseteq \Z^m$ be a positive affine seminormal monoid
such that $H^i_\mm(K[M])_{a} \neq 0$ for some $a\in\gp(M)$.
Then $a \in - \gp(M\cap F)\cap F$ for a face $F$ of $\cn(M)$ of dimension
$\leq i$.
In particular,
$$
H^i_\mm(K[M])_{a} =0\quad \text{if }a \not\in -\cn(M).
$$
\end{thm*}

As a consequence of this theorem and a careful analysis of the
groups $H^i_\mm(K[M])$ we obtain in \ref{seminice} that under the
hypothesis of the previous theorem $M$ is seminormal if and only if
$H^i_\mm(K[M])_a=0$  for all $i$ and all $a\in \gp(M)$ such that
$a\not\in -\cn(M)$. Note that this result has a variant for the
normalization, discussed in Remark \ref{remarknormality}.

Using further methods from commutative algebra we prove in Theorem
\ref{Hd_sn}:

\begin{thm*}
Let $M \subseteq \Z^m$ be a positive affine monoid of rank $d$
such that $K[M]$
satisfies $(S_2)$ and $H^d_\mm(K[M])_a=0$ for all $a\in
\gp(M)\setminus \bigl(-\cn(M)\bigr)$. Then $M$ is seminormal.
\end{thm*}

Since there are Cohen--Macaulay affine monoid rings which are not
seminormal (like $K[t^2,t^3]$), one can not omit the assumption
about the vanishing of the graded components  of $H^d_\mm(K[M])$
outside $-\cn(M)$.

In Section \ref{depthbound} we define the numbers
\begin{align*}
c_K(M) &= \sup\{i \in \Z: K[M\cap F] \text{ is Cohen--Macaulay for all
faces }
F,\ \dim F\leq i \},\\
 n(M) &= \sup\{i \in \Z: M\cap F \text{ is
normal for all faces } F,\ \dim F\leq i \}.
\end{align*}
By Hochster's theorem on normal monoids we have that $c_K(M) \geq
n(M)$. The main result \ref{maindepth1} in Section \ref{depthbound}
are the inequalities
$$
\depth R \geq c_K(M) \geq \min\{n(M)+1, d\}
$$
for an affine seminormal monoid $M \subseteq \Z^m$ of rank $d$.
Since seminormal monoids of rank $1$ are normal, we immediately get
that $\depth K[M] \geq 2$ if $d\ge 2$. In particular, $K[M]$ is
Cohen--Macaulay if $\rank M=2$.

One obtains a satisfactory result also for $\rank M=3$, which was
already shown in \cite{LI04} by different methods. In fact, in
Corollary \ref{cn_simp} we prove that $K[M]$ is Cohen--Macaulay
for a positive affine seminormal monoid $M \subseteq \Z^m$ with
$\rank M\leq 3$  if $K[M]$ satisfies $(S_2)$. One could hope that
$K[M]$ is always Cohen--Macaulay if $M$ is seminormal and $K[M]$
satisfies ($S_2$). But this is not the case and we present a
counterexample in \ref{exsemi}. The best possible result is given in
\ref{cn_simp}: ($S_2$) is sufficient for $K[M]$ to be Cohen-Macaulay
in the seminormal case if the cross-section of $\cn(M)$ is a simple
polytope.

In Section \ref{charp} we study the seminormality of affine monoid
rings with characteristic $p$ methods. The main observation is that
a positive affine monoid $M \subseteq \Z^m$ is seminormal if there
exists a field $K$ of characteristic $p$ such that $K[M]$ is
$F$-injective. This fact is a consequence of our analysis of the
local cohomology groups of positive affine monoid rings. In
\ref{primeprop} we give a precise description for which prime
numbers $p$ and fields of characteristic $p$, we have that $K[M]$ is
$F$-injective. Implicitly, this result was already observed by
Hochster and Roberts in \cite[Theorem 5.33]{HORO76}. In fact, if $M
\subseteq \Z^m$ is a positive affine seminormal monoid and $K$ is a
field of characteristic $p>0$, then the following statements are
equivalent:
\begin{enumerate}
\item
The prime ideal $(p)$ is not associated to the $\Z$-module $\gp(M)
\cap \R F/ \gp(M\cap F)$ for any face $F$ of $\cn(M)$;
\item
$R$ is $F$-split;
\item
$R$ is $F$-pure;
\item
$R$ is $F$-injective.
\end{enumerate}
As a direct consequence we obtain that $M$ is normal
if the equivalent statements hold for every field $K$ of characteristic $p>0$.

In the last section we present examples and counterexamples related
to the results of this paper. In particular, we will show that for
every simplicial complex $\Delta$ there exists a seminormal affine
monoid $M$ such that the only obstruction to the Cohen-Macaulay
property of $K[M]$ is exactly the simplicial homology of $\Delta$.
Choosing $\Delta$ as a triangulation of the real projective
plane we obtain an example whose Cohen-Macaulay property depends on
$K$. A similar result was proved by Trung and Hoa \cite{TRHO86}. Our
construction has the advantage of yielding a seminormal monoid $M$,
and is geometrically very transparent.

We are grateful to Aldo Conca for directing our attention toward the
results of Hochster and Roberts in \cite{HORO76}.
%------------------------------------------------------------------------
%
%
%
%------------------------------------------------------------------------
\section{Prerequisites}
We recall some facts from convex geometry. Let $X$ be a subset of
$\R^m$. The {\em convex hull} $\conv(X)$ of $X$ is the set of convex
combinations of elements of $X$. Similarly, the set $\cn(X)$ of
positive linear combinations of elements of $X$ is called the {\em
cone} generated by $X$. By convention $\cn(\emptyset)=\{0\}$ and
$\conv(\emptyset) = \emptyset$. A cone $C$ is called {\em positive}
(or {\em pointed}) if $0$ is the only invertible element in $C$. To an affine form
$\alpha$ on  $\R^m$ (i.~e.\ a polynomial of degree $1$) we associate
the affine hyperplane $H_{\alpha}= \alpha^{-1}(0)$, the closed
half-space $H_{\alpha}^+= \alpha^{-1}([0,\infty))$, and the open
half-space $H_{\alpha}^>= \alpha^{-1}(0,\infty)$. An intersection
$P=\bigcap_{i=1}^n H_{\alpha_i}^+$ of finitely many closed
half-spaces is called a {\em polyhedron}. A (proper) {\em face} of a
polyhedron $P$ is the (proper) intersection of $P$ with a hyperplane
$H_{\beta}$ such that $P \subseteq H_{\beta}^+$. Also $P$ is
considered as a face of itself. A {\em facet} is a maximal proper
face. Recall that there are only a finite number of faces. A {\em
polytope} is a bounded polyhedron. The set $\conv(F)$ is a polytope
for every finite subset $F$ of $\R^m$, and every polytope is of this
form. A {\em cone} is a finite intersection of half-spaces of the
form $H_{\alpha_i}^+$ where the $\alpha_i$ are linear forms
(i.\ e.\ homogeneous polynomials of degree $1$). The set $\cn(F)$ is a cone
for every finite subset $F$ of $\R^m$, and every cone is of this
form.
Let $P$ be a polyhedron and $F$ a face of $P$. Then we denote the
relative
interior of $F$ with respect to the subspace topology on the
affine hull of $F$ by $\relint F $. Note that $P$ decomposes into
the disjoint union $\relint F $ of the (relative) interiors of its
faces.
For more details on convex geometry we refer to the books of
Bruns and Gubeladze \cite{BRGUB},
Schrijver \cite{SCH98} and Ziegler \cite{ZI95}.

An {\em affine monoid} $M$ is a finitely generated commutative monoid
which can be embedded into $\Z^m$ for some $m \in \N$. We always
use $+$ for the monoid operation. In the literature $M$ is also
called an {\em affine semigroup} in this situation. We call $M$ {\em
positive} if $0$ is the only invertible element in $M$. Observe
that $M$ is positive if and only if $\cn(M)$ is pointed.

Let $K$ be a field and $K[M]$ be the $K$-vector space with
$K$-basis $X^a$, $a \in M$. The multiplication $X^a\cdot
X^b=X^{a+b}$ for $a,b \in M$ induces a ring structure on $K[M]$
and this $K$-algebra is called the affine monoid ring (or algebra)
associated to $M$. The embedding of $M$ into $\Z^m$ induces an
embedding of $K[M]$ into the Laurent polynomial ring $K[\Z^m]=
K[X_i^{\pm 1}:i=1,\dots,m]$ where $X_i$ corresponds to the $i$th
element of the canonical basis of $\Z^m$. Note that $K[M]$ is a
$\Z^m$-graded $K$-algebra with the property that $\dim_K
K[M]_a\leq 1$ for all $a \in \Z^m$. It is easy to determine the
$\Z^m$-graded prime ideals of $K[M]$. In fact every $\Z^m$-graded
prime ideal is of the form $\pp_F=(X^a: a\in M,\ a \not\in F)$ for
a unique face $F$ of $\cn(M)$ (see \cite[Theorem 6.1.7]{BRHE98}
for a proof). In particular, the prime ideals of height 1
correspond to the facets of $\cn(M)$.

Recall that a Noetherian domain $R$ is {\em normal} if it
is integrally closed in its field of fractions. The normalization
$\overline{R}$ of $R$ is the set of elements in the quotient field
of $R$ which are integral over $R$. An affine monoid $M$ is called
{\em normal}, if $z \in \gp(M)$ and $mz \in M$ for some $m\in \N$
imply $z \in M$. Here $\gp(M)$ is the group generated by $M$. It
is easy to see that $M$ is normal if and only if $M=\gp(M) \cap
\cn(M)$. If $M$ is an arbitrary submonoid of $\Z^m$, then its
{\em normalization} is the monoid $\overline{M}=\{z \in \gp(M) :
mz \in M \text{ for some } m \in \N \}$. By Gordan's lemma
$\overline M$ is affine for an affine monoid $M$.

Hochster \cite{HO72} proved that $K[M]$ is normal if and only if
$M$ is a normal monoid, in fact we have that
$\overline{K[M]}=K[\overline{M}]$. In particular,
Hochster showed that if
$M$ is normal, then $K[M]$ is a Cohen--Macaulay ring. One can
characterize the Cohen--Macaulay property of $K[M]$ in terms of
combinatorial and topological information associated to $M$. This
amounts to an analysis of the $\Z^m$-graded structure of the local
cohomology of $K[M]$; see Trung and Hoa \cite{TRHO86} for a
criterion of this type.

A Noetherian domain $R$ is called {\em seminormal} if
for an element $x$ in the quotient field $Q(R)$ of $R$ such that
$x^2,x^3 \in R$ we have $x \in R$. The {\em seminormalization}
$\+R$ of $R$ is the intersection of all seminormal subrings $S$
such that $R\subseteq S\subseteq Q(R)$.
An affine monoid $M$ is
called {\em seminormal}, if $z\in\gp(M)$, $z^2 \in M$ and $z^3 \in
M$ imply that $z \in M$. The seminormalization $\+M$ of $M$ is the
intersection of all seminormal monoids $N$ such that
$M\subseteq N \subseteq \gp(M)$. It can be shown that $\+M$ is again an affine
monoid. Hochster and Roberts \cite[Proposition 5.32]{HORO76}
proved that an affine monoid $M$ is seminormal if and only if
$K[M]$ is a seminormal ring. We frequently use the following
characterization of seminormal monoids.
See
Reid and Roberts \cite[Theorem 4.3]{RERO01}
for a proof for positive monoids, but this proof works also for arbitrary affine monoids.
\begin{thm}
\label{semicharacter} Let $M$ be an affine monoid $M \subseteq
\Z^m$. Then
$$
\+M=\bigcup_{F \text{ \rm face of } \cn(M) }  \gp(M\cap F) \cap
\relint F.
$$
In particular, $M$ is seminormal if and only if it equals the right
hand side of the equality.
\end{thm}

%------------------------------------------------------------------------
%
%
%
%------------------------------------------------------------------------
\section{Seminormality and Serre's condition $(S_2)$}
\label{s2section}
Let $R$ be a Noetherian ring and let $N$ be a finitely generated
$R$-module. Recall that $N$ satisfies {\em Serre's condition}
($S_k$) if
$$
\depth N_\pp \geq \min\{k,\dim N_\pp\}
$$
for all $\pp \in \Spec R$. Affine monoid rings trivially satisfy
($S_1$), since they are integral domains. We are interested in
characterizing $(S_2)$ for affine monoid rings.

While the validity of $(S_k)$ in $K[M]$ may depend on the field
$K$ for $k>2$, $(S_2)$ can be characterized solely in terms of
$M$, as was shown in \cite{SS90}.

Let $F_1, \dots, F_t$ be the facets
of $\cn(M)$ and let
$$
M_i=\{a \in \gp(M) : a+b \in M \text{ for some
} b \in M\cap F_i\}
$$
for $i=1,\dots,t$. Note that the elements of
$M_i$ correspond to the monomials
in the homogeneous localization $K[M]_{(\pp_{F_i})}$. We
set
$$
M'=\bigcap_{i=1}^t M_i.
$$

\begin{prop}\label{S2char}
\label{s2prime} Let $M$ be an affine monoid $M \subseteq \Z^m$, and $K$ a field. Then
the following statements are equivalent:
\begin{enumerate}
\item
$K[M]$ satisfies $(S_2)$;
\item
$M=M'$.
\end{enumerate}
\end{prop}

Observe that $(M')'=M'$. Thus $K[M']$ always satisfies $(S_2)$.
For seminormal monoids $M$ the equality $M=M'$ can be expressed in
terms of the lattices $\gp(M\cap F)$ as we will see in Corollary
\ref{s2cor}. First we describe the monoids $M_i$ under a slightly
weaker condition.

\begin{lem}[Proposition 4.2.6 in \cite{LI04}]
\label{lem427} Let $M\subseteq \Z^m$ be an affine monoid and
let $F_1,\dots,F_t$ be the facets of $\cn(M)$ with defining linear
forms $\alpha_1,\dots,\alpha_t$. Then:
\begin{enumerate}
\item
$ M_i \cap H_{\alpha_i}= \gp(M\cap F_i). $
\item
If $\gp(M)\cap \relint \cn(M) \subseteq M$, then
$$
M_i= \bigl(\gp(M) \cap H_{\alpha_i}^> \bigr) \cup \gp(M\cap F_i).
$$
\end{enumerate}
\end{lem}

\begin{proof}
Every element in $M_i$ is of the form $c=a-b$ for some $a\in M$ and
$b \in M\cap F_i$. Hence $\alpha_i(c)\geq 0$ with equality if and
only if $a \in M\cap F_i$. It follows that
$$
M_i \subseteq \bigl(\gp(M) \cap H_{\alpha_i}^>\bigr) \cup
\gp(M\cap F_i) \quad\text{and}\quad M_i \cap H_{\alpha_i}
\subseteq \gp(M\cap F_i).
$$
If $c \in \gp(M\cap F_i)$ and $c=a-b$ for some $a,b \in M\cap
F_i$, then clearly by the definition of $M_i$ we have that $c \in
M_i$. Thus we see that $M_i \cap H_{\alpha_i}= \gp(M\cap F_i)$.

For (ii) it remains to show that if
$c \in \gp(M) \cap H_{\alpha_i}^> $,
then $c \in M_i$. Pick $d \in M\cap \relint F_i$
such that $\alpha_j(c+d)>0$ for all $j \neq i$. Hence $c+d \in
\relint \cn(M)$. But
$$
c+d \in \gp(M)\cap \relint \cn(M) \subseteq M,
$$
by the additional assumption in (ii). Thus $c \in M_i$.
\end{proof}

In the following proposition we consider $\cn(M)$ as a face of itself.

\begin{prop}[Proposition 4.2.7 in \cite{LI04}]
\label{prop427} Let $M$ be an affine monoid $M \subseteq \Z^m$ and
let $F_1,\dots,F_t$ be the facets of $\cn(M)$. If $\gp(M)\cap \relint
\cn(M) \subseteq M$, then
$$
M'= \bigcup_{F \text{ \rm face of } \cn(M)} [\bigcap_{F\subseteq
F_i} \gp(M\cap F_i) \cap \relint F]
$$
with the convention that $\bigcap_{F\subseteq F_i} \gp(M\cap
F_i)=\gp(M)$ if $F=\cn(M)$.
\end{prop}
\begin{proof}
We apply \ref{lem427} several times. By assumption we have
$$
\gp(M) \cap \relint \cn(M) \subseteq M'.
$$
Let $F$ be a proper face of $\cn(M)$. Choose a facet $F_j$ with
defining linear form $\alpha_j$. Either $\relint F \subseteq F_j$
and thus
$$
\bigcap_{F \subseteq F_i} \gp(M\cap F_i) \cap \relint F \subseteq
\gp(M\cap F_j) \cap \relint F \subseteq M_j,
$$
or $\relint F$ is contained in $\cn(M)\cap H_{\alpha_j}^>$ and
$$
\bigcap_{F \subseteq F_i} \gp(M\cap F_i) \cap \relint F \subseteq
\gp(M)\cap H_{\alpha_j}^>\subseteq M_j.
$$
Hence
$$
\bigcap_{F \subseteq F_i} \gp(M\cap F_i) \cap \relint F \subseteq
M'.
$$
Note that
$$
M' \cap \relint \cn(M) \subseteq \gp(M) \cap \relint \cn(M).
$$
For a proper face $F$ of $\cn(M)$ it follows from \ref{lem427} that
$$
M' \cap \relint F \subseteq \bigcap_{F \subseteq F_i} M_i \cap
\relint F = \bigcap_{F \subseteq F_i} \gp(M\cap F_i) \cap \relint F.
$$
All in all we see that
\begin{equation*}
M'= \bigcup_{F \text{ \rm face of } \cn(M)} [\bigcap_{F \subseteq
F_i} \gp(M\cap F_i) \cap \relint F].\qedhere
\end{equation*}
\end{proof}

The equivalence of part (ii) and (iii) in the following corollary
was shown in Theorem 4.2.14 in \cite{LI04}.

\begin{cor}
\label{s2cor} Let $M \subseteq \Z^m$ be an affine monoid, $K$ a field, and let
$F_1,\dots,F_t$ be the facets of $\cn(M)$.  If $\gp(M)\cap \relint
\cn(M) \subseteq M$, then the following statements are equivalent:
\begin{enumerate}
\item
$K[M]$ satisfies $(S_2)$;
\item
$M=M'$;
\item
For all proper faces $F$ of $\cn(M)$ one has
$$
M \cap \relint F=\bigcap_{F \subseteq F_j} \gp(M\cap F_j)\cap
\relint F.
$$
\end{enumerate}
If $M$ is seminormal, then the following is equivalent to
\textup{(i) -- (iii)}:
\begin{enumerate}
\item[(iv)] $\gp(M\cap F)=\bigcap_{F \subseteq F_j} \gp(M\cap F_j)$.
\end{enumerate}
\end{cor}
\begin{proof}
The equivalence of (i) and (ii) was already stated in \ref{s2prime}.
The equivalence of (ii) and (iii) is an immediate consequence of
\ref{prop427}.

In the seminormal case one has
$\gp(F\cap M)\cap\relint F\subseteq M$
for all faces, so that (iv) implies (iii). For the converse
implication one uses the fact that $\gp(M\cap F)$ is generated by its
elements in $\relint(F)$ (see Bruns and Gubeladze \cite{BRGUB}).
\end{proof}

\begin{rem}
If $M$ is seminormal, then we know from \ref{semicharacter}, that
$\gp(M)\cap \relint \cn(M) \subseteq M$. Hence we can apply
\ref{lem427}, \ref{prop427} and \ref{s2cor} in this situation.
\end{rem}

The corollary shows that a seminormal monoid satisfies $(S_2)$ if
and only if the restriction of the groups $\gp(M\cap F)$ happens
only in the passage from $\cn(M)$ to its facets.

%-----------------------------------------------------------------------
%
%
%
%------------------------------------------------------------------------
\section{Local cohomology of monoid rings}
\label{localcohom} For the rest of the paper $K$ always denotes a
field, and $M\subseteq \Z^m$ is an affine positive monoid of rank
$d$. Recall that the seminormalization of $M$ is
$$
\+M=\bigcup_{F \text{ \rm face of } \cn(M) }  \gp(M\cap F) \cap
\relint F.
$$
and the normalization of $M$ is
$$
\overline{M}=\gp(M)\cap \cn(M).
$$
In this section we want to compute the local cohomology of $K[M]$
and compare it with the local cohomology of $K[\+M]$ and
$K[\overline{M}]$.

If $M$ is a positive affine monoid, then $K[M]$ is a $\Z^m$-graded
$K$-algebra with a unique graded maximal ideal $\mm$ generated by
all homogeneous elements of nonzero degree. By the local cohomology
of $K[M]$ we always mean the local cohomology groups
$H^i_\mm(K[M])$. Observe that $\+M$ and $\overline{M}$ are also
positive affine monoids. Since the $K$-algebras $K[\+M]$ and
$K[\overline M]$ are finitely generated modules over $K[M]$ and the
extensions of $\mm$ are primary to their maximal ideals, the local
cohomology groups of $K[\+M]$ and $K[\overline{M}]$ coincide with
$H^i_\mm(K[\+M])$ and $H^i_\mm(K[\overline{M}])$ respectively.
Because of this fact and to avoid cumbersome notation we always
write $H^i_\mm(K[\+M])$ and $H^i_\mm(K[\overline{M}])$ for the local
cohomology of $K[\+M]$ and $K[\overline{M}]$. The same applies to
$\Z^m$-graded residue class rings of $K[M]$.

In the following $R$ will always denote the ring $K[M]$, and thus
$\+R$ and $\overline{R}$ will stand for $K[\+M]$ and $K[\overline
M]$, respectively.

Let $F$ be a proper face of $\cn(M)$. Then
$\pp_F=(X^a: a\in M, a \not\in F)$ is a monomial prime ideal of $R$, and conversely, if $\pp$ is
a monomial prime ideal, then $F(\pp)=\R_+\{a\in M:X^a\notin
\pp\}$ is a proper face of $\cn(M)$. These two assignments set up
a bijective correspondence between the monomial prime ideals of
$R$ and the proper faces of $\cn(M)$.

Note that the natural embedding $K[M\cap F]\to K[M]$ is split by the
face projection $K[M]\to K[M\cap F]$ that sends all monomials in $F$
to themselves and all other monomials to $0$. Its kernel is $\pp_F$.
Therefore we have a natural isomorphism
$K[M]/\pp_F\cong K[M\cap F]$.

The next lemma states a crucial fact for the analysis of the
local cohomology of $R$. For this lemma and its proof we need
the following notation.
For $W \subseteq \Z^m$ we define
$$
-W=\{-a:a\in W\}.
$$
For a $\Z^m$-graded local Noetherian $K$-algebra $R$ with $R_0=K$
(like the monoid ring $K[M]$ for a positive affine monoid $M$)
and a $\Z^m$-graded $R$-module $N$ we set
$$
N^\vee=\Hom_K(N,K).
$$
(Here we mean by $\Hom_K(N,K)$ the homogenous homomorphisms from $N$ to $K$.)
Note that $N^\vee$ is again a $\Z^m$-graded $R$-module by setting
$$
(N^\vee)_a=\Hom_K(N_{-a},K) \text{ for } a \in \Z^m.
$$

\begin{lem}
\label{cmhelper}
Let $M \subseteq \Z^m$ be a positive affine monoid of rank $d$.
The $R$-module $\overline\omega$ of $R$ generated by the monomials
$X^b$ with $b \in \relint \cn(M) \cap \gp(M)$ is the canonical
module of the normalization $\overline{R}$. If $a \in \gp(M)$,
then
$$
H^i_\mm(\overline\omega)_{a} \cong
\begin{cases}
0 & \text{if } i<d \text{ or } a \not\in -\cn(M),\\
K & \text{if } i=d \text{ and } a \in    -\cn(M).
\end{cases}
$$
\end{lem}
\begin{proof}
Danilov and Stanley showed that $\overline\omega$ is the canonical
module of $\overline{R}$ (see Bruns and Herzog \cite[Theorem 6.3.5]{BRHE98} or
Stanley \cite{ST96}). Thus $\overline\omega$ is Cohen--Macaulay of
dimension $d$. This implies $H^i_\mm(\overline\omega)=0$ for
$i<d$. Furthermore, by graded local duality we have that
$$
H^d_\mm(\overline\omega)^\vee \cong
\Hom_{\overline{R}}(\overline\omega,\overline\omega) \cong
\overline{R}
$$
as $\Z^m$-graded modules. This concludes the proof.
\end{proof}

For the central proofs in this paper it is useful to extend the
correspondence between the faces of $\cn(M)$ and the monomial
prime ideals of $R$ to a bijection between the unions of faces of
$\cn(M)$ and the monomial radical ideals. If $\qq$ is a monomial
radical ideal, then we let $F(\qq)$ denote the union of the faces
$F(\pp)$ such that $\pp\supset \qq$, and if $F$ is the union of
faces, then the corresponding radical ideal $\qq_F$ of $R$ is
just the intersection of all monomial prime ideals $\pp_G$ such
that $G\subseteq F$.
We need the following lemma about monomial prime ideals of $R$.

\begin{lem}
\label{primehelper} Let $M \subseteq \Z^m$ be an affine monoid and
$F_1,\dots,F_t,G$ be faces of $\cn(M)$. Then:
\begin{enumerate}
\item
$\pp_{F_1\cap\dots\cap F_t}= \pp_{F_1} +\dots+ \pp_{F_t}$;
\item
$\pp_G+\bigcap_{i=1}^t \pp_{F_i}= \bigcap_{i=1}^t (\pp_G +\pp_{F_i})$.
\end{enumerate}
\end{lem}
\begin{proof}
Observe that the $\pp_{F_i}$ are $\Z^m$-graded ideals of $R$,
i.~e.\ they are monomial ideals in this ring. In other words,
their bases as $K$-vector spaces are subsets of the set of
monomials $X^a$, $a\in M$. Using this fact it is easy to check the
equalities claimed.
\end{proof}

We are ready to prove a vanishing result for the local cohomology
of seminormal monoid rings. Hochster
and Roberts \cite[Remark 5.34]{HORO76} already noticed  that certain ``positive'' graded
components of $H^i_\mm(R)$ vanish for a seminormal monoid. We can
prove a much more precise statement.

\begin{thm}\label{semihelper2}
Let $M \subseteq \Z^m$ be a positive affine seminormal monoid and
$R=K[M]$. If $H^i_\mm(R)_{a} \neq 0$ for $a\in\gp(M)$, then
$a \in - \overline{M\cap F}$ for a face $F$ of $\cn(M)$ of dimension
$\leq i$. In particular,
$$
H^i_\mm(R)_{a} =0\quad \text{if }a \not\in -\cn(M).
$$
\end{thm}
\begin{proof}
The assertion is trivial for $\rank M=0$. Thus assume that $d=\rank
M>0$. Since $M$ is seminormal, $\relint\cn(M) \cap
\gp(M)$ is contained in $M$. Thus $\overline\omega$, which as a
$K$-vector space is generated by the monomials $X^a$ with $a \in
\relint\cn(M) \cap \gp(M)$, is an ideal of $R$. Now
consider the exact sequence
$$
0 \to \overline\omega \to R \to R/ \overline\omega\to 0.
$$
By Lemma \ref{cmhelper} the long exact local cohomology sequences
splits into isomorphisms
\begin{eqnarray}
\label{eins} H^i_\mm(R) \cong H^i_\mm(R/ \overline\omega)
\quad\text{for } i<d-1
\end{eqnarray}
and the exact sequence
\begin{eqnarray}
\label{zwei} 0 \to H^{d-1}_\mm(R) \to H^{d-1}_\mm(R/
\overline\omega) \to H^{d}_\mm(\overline\omega) \to H^{d}_\mm(R) \to
0.
\end{eqnarray}
The local cohomology of $\overline\omega$ has been determined in
\ref{cmhelper}. Thus
$$
H^d_\mm(R)_{a}=0\quad \text{for } a \not\in -\cn(M).
$$
This takes care of the top local cohomology. For the lower
cohomologies we note that
$$
R/ \overline{\omega}\cong R/ \bigcap_{i=1}^t \pp_{F_i}
$$
where $F_1,\dots,F_t$ are the facets of $\cn(M)$.

Therefore it is enough to prove the following statement which
generalizes the theorem: \emph{let $\qq$ be a monomial radical ideal
of $R$; if $H^{i}_\mm(R/\qq)_{-a}\neq0$, then  $a \in
\overline{M\cap G}$ for a face $G\subseteq F(\qq)$ of dimension $\leq
i$}.

The case $\qq=(0)$ has already been reduced to the case
$\qq=\overline\omega$. So we can assume that $\qq\neq (0)$ and use
induction on $\rank M$ and on the number $t$ of minimal monomial
prime ideals $\pp_1,\dots,\pp_t$ of $\qq$.

If $t=1$, then $R/\pp_1\iso K[M\cap F(\pp_1)]$. Now we can apply
induction on $\rank M$.
Let $t>1$. We set $\qq'=\bigcap_{j=1}^{t-1} \pp_j$. Then we have the
standard exact sequence
$$
0\to R/\qq\to R/\qq'\oplus R/\pp_t \to R/(\qq'+\pp_t)\to 0.
$$
The local cohomologies of $R/\qq'$ and $R'=R/\pp_t$ are under
control by induction. But this applies to $R/(\qq'+\pp_t)$, too. In
fact, by Lemma \ref{primehelper} one has
$$
R/(\qq'+\pp_t)\cong R/\bigcap_{j=1}^{t-1}(\pp_t+\pp_j)\iso R'/\qq''
$$
where $\qq''=\bigcap_{j=1}^{t-1}\bigl((\pp_t+\pp_j)/\pp_t\bigr)$.
Thus $\qq''$ is a monomial radical ideal of $R'\cong K[M\cap
F(\pp_t)]$. Now it is enough to apply the long exact cohomology
sequence
\begin{equation*}
\dots \to H^{i-1}_\mm(R'/\qq'') \to H^{i}_\mm(R/\qq) \to
H^{i}_\mm(R/\qq') \oplus H^{i}_\mm(R')\to \cdots
\qedhere
\end{equation*}
\end{proof}

%------------------------------------------------------------------------------------
%------------------------------------------------------------------------------------
%------------------------------------------------------------------------------------
%------------------------------------------------------------------------------------
%------------------------------------------------------------------------------------
%------------------------------------------------------------------------------------
%------------------------------------------------------------------------------------
%------------------------------------------------------------------------------------

We describe a complex which computes the local cohomology of
$R$. Writing $R_F$ for the homogeneous localization
$R_{(\pp_F)}$, let
$$
L^{\pnt}(M)\colon 0 \to L^0(M) \to \dots \to L^t(M) \to \dots \to
L^d(M) \to 0
$$
be the complex with
$$
L^t(M)=\bigoplus_{F \text{ face of } \cn(M),\ \dim F=t} R_F
$$
and the differential $\partial \colon L^{t-1}(M) \to L^t(M)$
induced by
$$
\partial_{G,F} \colon R_G \to R_F
\text{ to be }
\begin{cases}
0 & \text{if } G \not\subset F,\\
\epsilon(G,F) \cdot\nat & \text{if } G \subset F,
\end{cases}
$$
where $\epsilon$ is a fixed incidence function on the face lattice
of $\cn(M)$ in the sense of \cite[Section 6.2]{BRHE98}. In
\cite[Theorem 6.2.5]{BRHE98} it was shown that $L^{\pnt}(M)$ is
indeed a complex and that for an $R$-module $N$ we have that
$$
H^i_\mm(N)=H^i(L^{\pnt}(M)\tensor_{R} N)\quad \text{ for all } i \geq
0.
$$

Next we construct another, ``smaller'' complex which will be
especially useful for the computation of the local cohomology of $R$
if $M$ is seminormal. Let
$$
\+L^{\pnt}(M)\colon 0 \to \+L^0(M) \to \dots \to \+L^t(M) \to \dots
\to \+L^d(M) \to 0
$$
be the complex with
$$
\+L^t(M)=\bigoplus_{F \text{ face of } \cn(M),\ \dim F=t}
K[-\overline{M\cap F}]
$$
and the differential $\+\partial \colon \+L^{t-1}(M) \to \+L^t(M)$
is induced by the same rule as $\partial$ above.

\begin{prop}
\label{semihelper1} Let $M \subseteq \Z^m$ be a positive affine
monoid, $a \in -\cn(M) \cap \gp(M)$. Then
$$
H^i_\mm(R)_{a} = H^i(\+L(M))_{a}
$$
\end{prop}
\begin{proof}
We know that $H^i_\mm(R)_{a}$ is the cohomology of the complex
$L(M)_{a}$. To prove the claim it suffices to determine
$(R_F)_{a}$ for a face $F$ of $\cn(M)$ if $a \in -\cn(M) \cap
\gp(M)$. It is an easy exercise to show that
$$
(R_F)_{a} = K[-\overline{M\cap F}]_{a}
$$
where one has to use the fact that $\overline{M\cap F}=\gp(M\cap
F)\cap F$.
\end{proof}
%------------------------------------------------------------------------------------
%------------------------------------------------------------------------------------
%------------------------------------------------------------------------------------
%------------------------------------------------------------------------------------
%------------------------------------------------------------------------------------
%------------------------------------------------------------------------------------
%------------------------------------------------------------------------------------
%------------------------------------------------------------------------------------

It follows from \ref{semihelper2} and \ref{semihelper1} that the
local cohomology of $\+R$ is a direct summand of the local
cohomology of $R$ as a $K$-vector space.

\begin{cor}
\label{semilocal} Let $M \subseteq \Z^m$ be a positive affine
monoid. Then
$$
\bigoplus_{a \in -\cn(M) \cap \gp(M)} H^i_\mm(R)_a \cong
\bigoplus_{a \in -\cn(M) \cap \gp(M)} H^i_\mm(\+R)_a =
H^i_\mm(\+R).
$$
\end{cor}

\begin{cor}
\label{semicor} Let $M \subseteq \Z^m$ be a positive affine monoid
of rank $d$. Then:
\begin{enumerate}
\item
If $R$ is Cohen--Macaulay, then $\+R$ is Cohen--Macaulay.
\item
If $\depth R \geq k$, then $\depth \+R \geq k$.
\item
If $R$ satisfies ($S_k$), then $\+R$ satisfies ($S_k$).
\end{enumerate}
\end{cor}
\begin{proof}
It is well-known that
the Cohen--Macaulay property and depth can be read off
the local cohomology groups.
This is also true for Serre's property ($S_k$)
since we have that
$R$
satisfies ($S_k$) if and only if
$\dim H^j_\mm(R)^\vee \leq j - k$
for
$j=0,\dots, \dim R -1$
and
an analogous characterization of  Serre's property ($S_k$)
for $\+R$.
(See Schenzel \cite{SCH} for a proof of the latter fact.)
\end{proof}

The results of this section allow us to give a cohomological
characterization of seminormality for positive monoid rings.

\begin{thm}
\label{seminice} Let $M \subseteq \Z^m$ be a positive affine monoid.
Then the following statements are equivalent:
\begin{enumerate}
\item
$M$ is seminormal;
\item
$H^i_\mm(R)_a=0$  for all $i$ and all $a\in \gp(M)$ such that
$a\not\in -\cn(M)$.
\end{enumerate}
\end{thm}
\begin{proof}
Consider the sequence
$$
0 \to R \to \+R \to \+R/R \to 0
$$
of finitely generated $\Z^m$-graded $R$-modules. Observe that
$H^i_\mm(R)_a\cong H^i_\mm(\+R)_a$ for $a \in -\cn(M)$. Thus it
follows from the long exact cohomology sequence
$$
\dots \to H^i_\mm(R) \to H^i_\mm(\+R) \to H^i_\mm(\+R/R) \to
\dots
$$
that $H^i_\mm(R)_a=0$ for $a\not\in -\cn(M)$ and all $i$ if and only
if $H^i_\mm(\+R/R)_a=0$ for all $a $ and all $i$. This is equivalent
to $\+R/R=0$. Hence $R$ and, thus, $M$ are seminormal.
\end{proof}

\begin{rem}
\label{remarknormality}
The previous results have variants for the normalization. If we
restrict the direct sum in Corollary \ref{semilocal} to those $a$
that belong to $-\relint\cn(M)\cap\gp(M)$ then the local
cohomology of $K[\+M]$ must be replaced by that of $\overline{R}$.
Moreover, the local cohomology of $R$ vanishes in all degrees $a$
outside $-\relint\cn(M)$ if and only if $M$ is normal.

This follows by completely analogous arguments since
we have that
$H_\mm^i(\overline{R})=0$ for $i<d$ and
$H_\mm^d(\overline{R})_a\neq 0$ if and only if
$a\in-\relint\cn(M)\cap\gp(M)$.
\end{rem}
%------------------------------------------------------------------------------------
%------------------------------------------------------------------------------------
%------------------------------------------------------------------------------------
%------------------------------------------------------------------------------------
%------------------------------------------------------------------------------------
%------------------------------------------------------------------------------------

With different methods than those used so far we can prove another
seminormality criterion. It involves only the top local cohomology
group, but needs a stronger hypothesis on $M$.

\begin{thm}\label{Hd_sn}
Let $M \subseteq \Z^m$ be a positive affine monoid of rank $d$. If
$R$ satisfies $(S_2)$ and $H^d_\mm(R)_a=0$ for all $a\in
\gp(M)\setminus \bigl(-\cn(M)\bigr)$, then $M$ is seminormal.
\end{thm}
\begin{proof}
The assumption and \ref{semilocal} imply that the $d$th local
cohomology of $R$ and $\+R$ coincide as $R$-modules. Since $M$ and
therefore $\+M$ are positive, there exists a $\Z$-grading on $R$ and
$\+R$ such that both $K$-algebras are generated in positive degrees.
We choose a common Noether normalization $S$ of $R$ and $\+R$ with
respect to this $\Z$-grading.

Since $R$ satisfies $(S_2)$ it is a reflexive $S$-module. By
\ref{semicor} also $\+R$ satisfies $(S_2)$ and is a reflexive
$S$-module. In the following let $\omega_S$ be the canonical module
of $S$ which is in our situation just a shifted copy of $S$ with
respect to the $\Z$-grading. By graded local duality and reflexivity
we get the following chain of isomorphisms of graded $S$-modules:
\begin{align*}
R & \cong  \Hom_S(\Hom_S(R,S),S) \cong  \Hom_S(\Hom_S(R,\omega_S),\omega_S)\\
& \cong \Hom_S(H^d_\mm(R)^\vee,\omega_S) \cong  \Hom_S(H^d_\mm(\+R)^\vee,\omega_S)\\
& \cong \Hom_S(\Hom_S(\+R,\omega_S),\omega_S) \cong  \Hom_S(\Hom_S(\+R,S),S)\\
& \cong \+R.
\end{align*}
Hence $M=\+M$ is seminormal.
\end{proof}

Again we can obtain a similar normality criterion if we replace
$-\cn(M)$ by $-\relint\cn(M)$ in Theorem \ref{Hd_sn}. In the rest of
this section we further analyze the local cohomologies of~$R$.

\begin{prop}\label{HdHd-1}
Let $M$ be seminormal and $R=K[M]$. Then:
\begin{enumerate}
\item $H_\mm^d(R)_{-a}\neq 0$ (and so of $K$-dimension $1$)
$\Leftrightarrow$ $a\in\overline M\setminus\bigcup_F \gp(M\cap F)$ where $F$
runs through the facets of $\cn(M)$.
\item $H_\mm^{d-1}(R)_{-a}\neq 0$ $\Leftrightarrow$ $a\in
\partial\cn(M)\cap \bigcup_F \gp(M\cap F)$ and
$\dim_K H_\mm^{d-1}(R/\overline\omega)_{-a}\ge 2$ where $F$ again
runs through the facets of $\cn(M)$.
\item $R$ is Cohen--Macaulay $\Leftrightarrow$ $R/\overline\omega$ is
Cohen--Macaulay and $\dim_K H_\mm^{d-1}(R/\overline\omega)_{-a}\le 1$
for all $a\in\gp(M)$.
\end{enumerate}
\end{prop}

\begin{proof}
$H_\mm^d(R)$ is the cokernel of the map
$$
\bigoplus_{F \text{ facet of } \cn(M)} K[-\overline{M\cap F}]\ \to\
K[-\cn(M)\cap\gp(M)],
$$
which implies (i).

(ii) We have the exact sequence \eqref{zwei}
$$
0 \to H^{d-1}_\mm(R) \to H^{d-1}_\mm(R/ \overline\omega) \to
H^{d}_\mm(\overline\omega) \xrightarrow{\pi} H^{d}_\mm(R) \to 0.
$$
Lemma \ref{cmhelper} and (i) yield the $\Z^m$-graded structure of
$\Ker\pi$: its nonzero graded components have dimension $1$ and live
in exactly the degrees $-a$ with $a\in \cn(M)\cap\bigcup_F \gp(M\cap
F)$.

On the other hand, $H_\mm^{d-1}(R/\overline\omega)$ can have
non-zero components only in these degrees, as follows from the
generalization of Theorem \ref{semihelper2} stated in its proof.
Thus $H_\mm^{d-1}(R)$ is limited to these degrees and is non-zero at
$-a$ if and only if $\dim_K H_\mm^{d-1}(R/\overline\omega)_{-a}\ge
2$.

(iii) The isomorphisms $H_\mm^i(R)\cong H_\mm^i(R/\overline\omega)$
for $i<d-1$ reduce the claim immediately to (ii).
\end{proof}

Having computed the $d$-th local cohomology of $K[M]$, we can easily
describe the Gorenstein property of $K[M]$ in combinatorial terms:

\begin{cor}
Let $M$ be seminormal and $R=K[M]$ Cohen--Macaulay. For each facet
$F$ of $\cn(M)$ let $\gamma_F$ denote the index of the group
extension $\gp(M\cap F)\subset \gp(M)\cap \R F$, and $\sigma_F$ the
unique $\Z$-linear form on $\gp(M)$ such that $\sigma_F(\overline
M)=\Z_+$ and $\sigma_F(x)=0$ for all $x\in F$. Then the following
are equivalent:
\begin{enumerate}
\item $R$ is Gorenstein;
\item
\begin{enumerate}
\item $\gamma_F\le 2$ for all facets $F$ of $\cn(M)$;
\item there exists $b\in \overline M$ such that
\begin{alignat*}2
&b\in F\setminus \gp(M\cap F),\quad &&\text{if }\gamma_F=2,\\
&\sigma_F(b)=1,&&\text{else}.
\end{alignat*}
\end{enumerate}
\end{enumerate}
\end{cor}

\begin{proof}
The multigraded support of the $K$-dual $\omega_R$ of $H_\mm^d(R)$
is $N=\overline M\setminus \bigcup_F \gp(M\cap F)$, and $R$ is
Gorenstein if and only if there exists $b\in \overline M$ such that $N=b+M$,
or, equivalently, $\omega_R=RX^b$. It remains to be shown that such
$b$ exists if and only if the conditions in (ii) are satisfied. We
leave the exact verification to the reader. (Note that for each
facet $F$ there exists $c\in \overline M\cap \relint\cn(M)\subset M$
such that $\sigma_F(x)=1$.)
\end{proof}
%------------------------------------------------------------------------------------
%------------------------------------------------------------------------------------
%------------------------------------------------------------------------------------
%------------------------------------------------------------------------------------
%------------------------------------------------------------------------------------
%------------------------------------------------------------------------------------

Finally, we give an interpretation of $H^i_\mm(R)_{-a}$ for $a \in
\cn(M)$ which will be useful in later sections.

\begin{rem}
\label{topological} Let $M \subseteq \Z^m$ be a positive affine
monoid of rank $d$, $i \in \{0,\dots, d\}$, $a \in \cn(M)\cap
\gp(M)$ and $\Fc(M,a)=\{F \text{ face of } \cn(M) : a \in
\overline{M\cap F}\}$. Then $H^i_\mm(R)_{-a}$ is the $i$th
cohomology of the complex
$$
\mathcal{C}^{\pnt}(M,a)\colon 0 \to \mathcal{C}^{0}(M,a) \to \dots
\to \mathcal{C}^{t}(M,a) \to \dots \to \mathcal{C}^{d}(M,a) \to 0,
$$
where
$$
\mathcal{C}^{t}(M,a)=\bigoplus_{G \in \Fc(M,a) ,\ \dim G=t} Ke_G
$$
and the differential is given by $\partial(e_G)=\sum_{F \text{ face
of } \cn(M),\ G \subseteq F} \epsilon(G,F)e_F$. (Here $\epsilon$ is
the incidence function on the face lattice of $\cn(M)$ which we
fixed above to define the complex $L^{\pnt}(M)$.)
\end{rem}

\begin{thm}
\label{vanishpiece} Let $M \subseteq \Z^m$ be a positive affine
monoid of rank $d$ and $a \in \cn(M)\cap\gp(M)$. If the set
$\Fc(M,a)=\{F \text{ face of } \cn(M) : a \in \overline{M\cap F}\}$
has a unique minimal element $G$, then
$$
H^i_\mm(R)_{-a}=0\quad \text{for all } i=0,\dots, d-1.
$$
\end{thm}

\begin{proof}
It follows from \ref{topological} that $H^i_\mm(R)_{-a}$ is the
$i$th cohomology of the complex
$$
\mathcal{C}^{\pnt}(M,a)\colon 0 \to \mathcal{C}^{0}(M,a) \to \dots
\to \mathcal{C}^{t}(M,a) \to \dots \to \mathcal{C}^{d}(M,a) \to 0,
$$
where
$$
\mathcal{C}^{t}(M,a)=\bigoplus_{G \in \Fc(M,a) ,\ \dim G=t} Ke_G.
$$
By taking a cross section of $\cn(M)$ we see that the face lattice
of $\cn(M)$ is the face lattice of a polytope (see
\cite[Proposition 6.1.8]{BRHE98}). If $\Fc(M,a)$ has a unique
minimal element, then this set is again the face lattice of a
polytope $P$, as can be seen from Ziegler \cite[Theorem 2.7]{ZI95}. Note
that if $\Fc(M,a)$ has only one element, then $P$ is the empty set.
But this can only happen if $\Fc(M,a)=\{\cn(M)\}$ and then we have
homology only in cohomological degree $d$. If $\Fc(M,a)$ has more
than one element, then $\mathcal{C}^{\pnt}(M,a)$ is the $K$-dual
of a cellular resolution which computes the singular cohomology of
$P$. A nonempty polytope is homeomorphic to a ball and thus the
complex $\mathcal{C}^{d}(M,a)$ is acyclic. Hence in this case
$H^i_\mm(R)_{-a}=0$ for $i=0,\dots,d$, and this concludes the
proof.
\end{proof}

The following corollary collects two immediate consequences of
Theorem \ref{vanishpiece}.

\begin{cor}\label{van_cor}
Let $M$ be a positive affine monoid.
\begin{enumerate}
\item Let $F$ be the unique face of $\cn(M)$ such that $a \in
\relint F$. If $a \in \gp(M\cap F)$, then
$$
H^i_\mm(R)_{-a}=0\quad \text{ for all } i=0,\dots, d-1.
$$
\item Suppose that $M$ is seminormal. Then $R$ is Cohen--Macaulay if
$\Fc(M,a)$ has a unique minimal element for all $a\in \overline{M}$.
\end{enumerate}
\end{cor}

Finally we note that nonzero lower local cohomologies must be large
in the seminormal case.

\begin{prop}
Let $M \subseteq \Z^m$ be a positive affine seminormal monoid. If
$H^i_\mm(R)_a \neq 0$ for some $a \in - \cn(M)$, then $\dim_K
H^i_\mm(R)=\infty$. In particular, if a seminormal monoid is
Buchsbaum, then it must be Cohen--Macaulay.
\end{prop}

\begin{proof}
If $H^i_\mm(R)_a \neq 0$, then the complex $\mathcal{C}^{\pnt}(M,a)$
has nontrivial cohomology in degree $i$. Consider the multiples $ka$
for $k \in \N$. If $a \in \overline{M\cap F}=\gp(M\cap F)\cap F$,
then $ka \in \overline{M\cap F}$ for all $k \in \N$. If $a \not\in
\overline{M\cap F}$, then there exist infinitely many $k$ such that
$ka \not\in \overline{M\cap F}$. Since the face lattice of $\cn(M)$
is finite we can choose a sequence $(k_n)_{n\in \N}$ such that
$k_n<k_{n+1}$ and $a \not\in \overline{M\cap F}$ if and only if $k_n
a \not\in \overline{M\cap F}$. Thus
$\mathcal{C}^{\pnt}(M,a)=\mathcal{C}^{t}(M,k_n a)$ for all $n\geq 0$
which implies $H^i_\mm(R)_{k_na} \neq 0$. Hence $\dim_K
H^i_\mm(R)=\infty$.

If $R$ is Buchsbaum, then $\dim_K H^i_\mm(R)<\infty$ for all
$i<d$. Thus the local cohomology must vanish in this case for
$i<d$ which implies that $R$ is already Cohen--Macaulay.
\end{proof}

\section{The Cohen--Macaulay property and depth}
\label{depthbound}

If a seminormal monoid $M$ fails to be normal by the smallest
possible margin, then $K[M]$ is Cohen--Macaulay
as the following result shows:

\begin{prop}\label{normalfacet}
Let $M$ be seminormal such that $M\cap F$ is normal for each facet
$F$ of $\cn(M)$. Then $R$ is Cohen--Macaulay.
\end{prop}

\begin{proof}
It is enough to show that $\Fc(M,a)$ has a unique minimal element
for all $a\in \overline{M}$. Let $a\in \gp(M)\cap \cn(M)$. If
$a\notin \overline{M\cap G}$ for all facets $G$ of $\cn(M)$, then
$\cn(M)$ is the unique minimal element of $\Fc(M,a)$.

Otherwise we have $a\in \overline{M\cap G}=M\cap G\subset M$ for
some facet $G$ of $\cn(M)$. We choose the unique face $F'$ of
$\cn(M)$ with $a\in \relint F'$. It follows that $a\in F'\cap M$,
and $F'$ is the unique minimal element of $\Fc(M,a)$.
\end{proof}

\begin{rem}
Another, albeit more complicated proof of the proposition can be
given as follows. The main result of Brun, Bruns and R\"omer \cite{BBR05} implies for
$R=K[M]$ that
\begin{enumerate}
\item $R/\overline\omega$ is Cohen--Macaulay,
\item $H_\mm^d(R/\overline\omega)=\bigoplus_F H_\mm^{\dim F}(K[M\cap
F])$ where $F$ runs through the proper faces of $\cn(M)$,
\end{enumerate}
provided all the rings $K[M\cap F]$ are Cohen--Macaulay. If they are
even normal, then the local cohomology modules in (ii) do not
``overlap'' because the $\Z^m$-graded support of $H_\mm^{\dim
F}(K[M\cap F])$ is restricted to $-\relint F$, and the relative
interiors of faces are pairwise disjoint. Now we can conclude from
Proposition \ref{HdHd-1} that $R$ is Cohen--Macaulay.

In general, without normality of the facets the local cohomology
modules in (ii) will overlap (see Example \ref{exsemi}). This limits
all attempts to prove stronger assertions about the Cohen--Macaulay
property in the seminormal case.
\end{rem}

Using the results and techniques of Section \ref{localcohom}, we
can give lower bounds for the depth of seminormal monoid rings.
Let $M \subseteq \Z^m$ be an affine seminormal monoid. We define
\begin{align*}
c_K(M) &= \sup\{i \in \Z: K[M\cap F] \text{ is Cohen--Macaulay for all
faces }
F,\ \dim F\leq i \},\\
 n(M) &= \sup\{i \in \Z: M\cap F \text{ is
normal for all faces } F,\ \dim F\leq i \}.
\end{align*}
Observe that if $M\cap F$ is normal for a face $F$ of $\cn(M)$, then
also $M\cap G$ is normal for all faces $G\subseteq F$ of $\cn(M)$.
Hence it would be enough to consider all $i$-dimensional faces of
$\cn(M)$ in the definition of $n(M)$. However, as we will see in
Section \ref{ExandCEx}, this is not true for the Cohen--Macaulay
property.

\begin{thm}
\label{maindepth1} Let $M \subseteq \Z^m$ be an affine seminormal
monoid of rank $d$, and $R=K[M]$. Then
$$
\depth R \geq c_K(M) \geq \min\{n(M)+1, d\}.
$$
\end{thm}

\begin{proof}
The proof of the first inequality follows essentially the same idea
as that of Theorem \ref{semihelper2}.

The assertion is trivial for $\rank M=0$. Thus assume that $d=\rank
M>0$. There is nothing to prove if $c_K(M)=d$. So we can assume that
$c_K(M)<d$. Since $M$ is seminormal, we can again use the exact
sequence
$$
0 \to \overline\omega \to R \to R/ \overline\omega\to 0.
$$
Since $\depth \overline\omega= d$ according to Lemma \ref{cmhelper},
it is enough to show that $\depth R/\overline\omega\ge c_K(M)$.

Again we write $\overline\omega=\bigcap_{j=1}^t\pp_{F_j}$ where
$F_1,\dots,F_t$ are the facets of $\cn(M)$. However, contrary to
Theorem \ref{semihelper2}, the bound does not hold for arbitrary
residue class rings with respect to monomial radical ideals $\qq$,
since the combinatorial structure of the set $F(\qq)$ may contain
obstructions.

Therefore we order the facets $F_1,\dots,F_t$ in such a way that
they form a shelling sequence for the face lattice of $\cn(M)$.
Such a sequence exists by the Brugesser-Mani theorem (applied to a
cross section polytope of $\cn(M)$). See \cite[Lecture 8]{ZI95}.
The generalization of the first inequality of the theorem to be
proved is the following: \emph{let $F_1,\dots,F_t$ be a shelling
sequence for $\cn(M)$ and let $u\in\{1,\dots,t\}$, then $\depth
R/\qq\ge \min\{d-1,c_K(M)\}$ for $\qq=\bigcap_{j=1}^u \pp_{F_j}$.}

If $u=1$, then $R/\pp_1\iso K[M\cap F(\pp_1)]$. Now we can apply
induction on $\rank M$.
Let $u>1$. We set $\qq'=\bigcap_{j=1}^{u-1} \pp_{F_j}$. Again we
have the standard exact sequence
$$
0\to R/\qq\to R/\qq'\oplus R/\pp_{F_u}\to R/(\qq'+\pp_{F_u})\to
0.
$$
Therefore
$$
\depth R/\qq \ge \min\{ 1+\depth(R/(\qq'+\pp_{F_u}),\depth
R/\pp_{F_u},\depth R/\qq'\}.
$$
By induction on $u$ we have $\depth R/\pp_{F_u},\depth R/\qq'
\ge\min\{d-1,c_K(M)\}$.

Now the crucial point is that
$F_u \cap \bigcup_{j=1}^{u-1}F_j= \bigcup_{j=1}^{u-1}F_u\cap F_j$ is the union of certain facets
$G_1,\dots,G_v$ of $F_u$ that form the starting segment of a
shelling sequence for $F_u$ (by the very definition of a shelling).
As in the proof of Theorem \ref{semihelper2} we have
$$
R/(\qq'+\pp_{F_u})\cong
R/\bigcap_{j=1}^{t-1}(\pp_{F_u}+\pp_{F_j})\iso R'/\qq''
$$
where
$\qq''=\bigcap_{j=1}^{t-1}\bigl((\pp_{F_u}+\pp_{F_j})/\pp_{F_u}\bigr)$.
Therefore $\qq''$ is the radical ideal of $R'=K[M\cap F_1]$
corresponding to the union $G_1,\dots,G_v$. By induction we have
$$
\depth R'/\qq''\ge \min\{ d-2,c(M\cap F_u)\}
\ge \min\{ d-2,c_K(M)\},
$$ and this completes
the proof for the inequality $\depth R\ge c_K(M)$.

By Hochster's theorem the second inequality holds if $M$ itself is
normal. Suppose that $n(M)<d$ and let $F$ be a face of dimension
$n(M)+1$. Then we must show that $K[M\cap F]$ is Cohen--Macaulay.
Thus the second inequality reduces to the claim that $R$ is
Cohen--Macaulay if the intersections $F\cap M$ are normal for all
facets $F$ of $\cn(M)$ (and $M$ is seminormal). This has been shown
in Proposition \ref{normalfacet}.
\end{proof}

There is a general lower bound for the depth of seminormal monoid
rings of rank $\geq 2$. It follows from the proposition since
seminormal monoids of rank $1$ are normal.

\begin{cor}
Let $M \subseteq \Z^m$ be an affine seminormal monoid of rank $d\geq
2$. Then
$$
\depth R \geq 2.
$$
In particular, $R$ is Cohen--Macaulay if $d=2$.
\end{cor}

One could hope that seminormality plus some additional assumptions
on $M$ already imply the Cohen--Macaulay property of $R$. But most
time this is not the case as will be discussed in Example
\ref{exsemi}. However, we will now show that Serre's condition
$(S_2)$ implies the Cohen--Macaulay property of $R$ if $\cn(M)$ is
a simple cone (to be explained below).
More generally, we want to show that simple faces of
$\cn(M)$ cannot contain an obstruction to the Cohen--Macaulay
property in the presence of $(S_2)$.

Let $F$ be a proper face of $\cn(M)$. We call the face $F$ {\em
simple} if the partially ordered set $\{G \text{ face of } \cn(M):
F\subseteq G\}$ is the face lattice of a simplex. Observe that by
\cite[Theorem 2.7]{ZI95} the latter set is always the face lattice
of a polytope, because the face lattice of $\cn(M)$ is the face
poset of a cross section of $\cn(M)$.
Let $F$ be a simple face of $\cn(M)$.
It is easy to see that
every face $G$ of $\cn(M)$ containing the simple face $F$ is also simple.

We say that $\cn(M)$ is simple if a cross section polytope of
$\cn(M)$ is a simple polytope.
This amounts to the simplicity of all the edges
of $\cn(M)$. (Note that the apex $\{0\}$ is a simple face if and
only if $\cn(M)$ is a simplicial cone.)

\begin{prop}
\label{simplecrit} Let $M \subseteq \Z^m$ be a positive affine
seminormal monoid such that $R$ satisfies $(S_2)$. Let $a \in
\gp(M)\cap \cn(M)$ and $a \in \relint F$ for a proper face $F$ of
$\cn(M)$. If $H^i_\mm(R)_{-a} \neq 0$ for some $i$, $0\leq i \leq
d-1$, then $F$ is not a simple face of $\cn(M)$.
\end{prop}
\begin{proof}
Assume that $F$ is a simple face. Consider the intersection
$$
H = \bigcap_{G \text{ face of } \cn(M),\ F\subseteq G,\ a \in
\overline{M\cap G}} G
$$
which is a simple face containing $F$ because $F$ is simple. Let
$F_1,\dots,F_t$ be the facets of $\cn( M)$. For each facet $F_j$
such that $H\subseteq F_j$ there exists a face $G$ of $\cn(M)$
with $F\subseteq G$, $a \in \overline{M\cap G}$ such that
$G\subseteq F_j$ because $H$ is simple. This follows from the fact
that the partially ordered set $\{L: L \text{ is a face of }
\cn(M)$, $H\subseteq L \}$ is the face poset of a simplex,
and for a simplex the claim is trivially true. We observe
that $a \in \gp(M\cap G) \subseteq \gp(M\cap F_j)$ for those
facets $F_j$ with $H \subseteq F_j$.
By Corollary \ref{s2cor} we have
$$
\gp(M\cap H)=\bigcap_{H\subseteq F_j} \gp(M\cap F_j).
$$
Therefore $a \in \gp(M\cap H)\cap H = \overline{M\cap H}$.

All in all we get that the set $\Fc(M,a)=\{L \text{ face of } \cn(M) :
a \in \overline{M\cap L}\}$ has the unique minimal element $H$, and
\ref{vanishpiece} implies that
$$
H^i_\mm(R)_{-a} = 0
$$
which is a contradiction to our assumption. Thus $F$ is not a simple
face of $\cn(M)$.
\end{proof}

The latter result gives a nice Cohen--Macaulay criterion in terms of
$\cn(M)$ for a seminormal monoid. It implies Theorem 4.4.7 in
\cite{LI04}, and can be viewed as a variant of the theorem by Goto,
Suzuki and Watanabe \cite{GSW76} by which $(S_2)$ implies the
Cohen--Macaulay property of $R$ if $\cn(M)$ is simplicial.

\begin{cor}\label{cn_simp}
Let $M \subseteq \Z^m$ be a positive affine seminormal monoid such
that $R=K[M]$ satisfies $(S_2)$ and such that $\cn(M)$ is a simple
cone. Then $R$ is Cohen--Macaulay for every field $K$. In
particular, if $\rank M\leq 3$, then $R$ is Cohen--Macaulay.
\end{cor}
\begin{proof}
Every proper face of $\cn(M)$, with the potential exception of
$\{0\}$, is simple. Thus it follows from \ref{semihelper2} and
\ref{simplecrit} that $H^i_\mm(R)_{-a}=0$ for $a\neq0$ and
$i=0,\dots,d-1$. For $a=0$ this results from Corollary
\ref{van_cor}. Hence $R$ is Cohen--Macaulay.

If $\rank M \leq 3$, then the cross section of $\cn(M)$ is a
polytope of dimension $\leq 2$, which is always simple. Thus we can
apply the corollary.
\end{proof}

We will point out in Remark \ref{best_poss} that the corollary is
the best possible result if one wants to conclude the
Cohen--Macaulay property of $R$ only from the seminormality of $M$
and the validity of ($S_2$).

%------------------------------------------------------------------------
%
%
%
%------------------------------------------------------------------------
\section{Seminormality in characteristic $p$}
\label{charp}
In this section we study local cohomology properties of seminormal
rings in characteristic $p>0$. Let $K$ be a field with $\chara
K=p>0$. In this situation we have the Frobenius homomorphism
$F\colon R \to R,\ f \mapsto f^p$. Through this homomorphism $R$ is
an $F(R)$-module. Now $R$ is called {\em $F$-injective} if the
induced map on the local cohomology $H^i_{\mm}(R)$ is injective for
all $i$. It is called {\em $F$-pure} if the extension $F(R) \to R$
is pure, and {\em $F$-split} if $F(R)$ is a direct $F(R)$-summand of
$R$. In general we have the implications
$$
\text{$F$-split } \implies \text{ $F$-pure } \implies \text{
$F$-injective}.
$$
For example see \cite{BRHE98} for general properties of these
notions.

\begin{prop}
\label{charp1}
Let $M \subseteq \Z^m$ be a positive affine monoid. If there exists
a field $K$ of characteristic $p$ such that $R$ is $F$-injective,
then $M$ is seminormal.
\end{prop}

\begin{proof}
Assume that there exists an $a \in \gp(M)$, $a \not\in -\cn(M)$, and
an $i\in \{0,\dots, \rank M\}$ such that $H^i_\mm(R)_a\neq 0$. Since
$R$ is $F$-injective, it follows that $H^i_\mm(R)_{p^m \cdot a}\neq
0$ for all $m\in \N$. Write $R=S/I_M$ as a $\Z^m$-graded quotient of
a polynomial ring $S$. Then by graded local duality $H^i_\mm(R)^\vee
\cong \Ext_S^{n-i}(R,\omega_S)$ is a finitely generated
$\Z^m$-graded $R$-module. This implies that $H^i_\mm(R)_{p^m \cdot
a}= 0$ for $m \gg 0$, which is a contradiction. Thus $H^i_\mm(R)_{a}=
0$ for all $a \not\in \cn(M)$. It follows from \ref{seminice}
that $M$ is seminormal.
\end{proof}

If $M$ is seminormal there exist only finitely many prime numbers
such that $R$ is not $F$-injective. Moreover, we can characterize
this prime numbers precisely.

\begin{prop}
\label{primeprop} Let $M \subseteq \Z^m$ be a positive affine
seminormal monoid and let $K$ be a field of characteristic $p>0$.
Then the following statements are equivalent:
\begin{enumerate}
\item
The prime ideal $(p)$ is not associated to the $\Z$-module
$\gp(M) \cap \R F/ \gp(M\cap F)$ for any face $F$ of $\cn(M)$;
\item
$R$ is $F$-split;
\item
$R$ is $F$-pure;
\item
$R$ is $F$-injective.
\end{enumerate}
\end{prop}

\begin{proof}
(i) $\Rightarrow$ (ii) We show by a direct computation that
$F(R)=K^p[pM]$ is a direct $K^p[pM]$-summand of $R$. Since $K^p[pM]$
is a $K^p[pM]$-summand of $K[pM]$, it is enough to show that $K[pM]$
is a direct $K[pM]$-summand of $R$. The monoid $M$ is the disjoint
union of the residue classes modulo $p\gp(M)$. This induces a direct
sum decomposition of $R$ as a $K[pM]$-module. We claim that $pM$ is
the intersection of $M$ and $p\gp(M)$. This will show the remaining
assertion.

To prove the claim we have only to show that an element $w$ of the
intersection of $M$ and $p \gp(M)$ is an element of $pM$. Write
$w=pz$ for some $z \in \gp(M)$. We have that $w \in \relint F$ for a
face $F$ of $\cn(M)$. Then $p$ annihilates the element $z \in \gp(M)
\cap \R F$ modulo $\gp(M\cap F)$. By assumption $p$ is a
nonzero-divisor on that module. Thus $z \in \gp(M\cap F)$. Since $z
\in \gp(M\cap F) \cap \relint F$ and $M$ is seminormal we have that
$z\in M$ by \ref{semicharacter}. Hence $w \in pM$ as desired.

(ii) $\Rightarrow$ (iii)
and
(iii) $\Rightarrow$ (iv)
: This holds in general as was remarked
above.

(iv) $\Rightarrow$ (i) Assume that $(p)$ is associated to some of
the $\Z$-modules $\gp(M) \cap \R F/ \gp(M\cap F)$ for the faces $F$ of
$\cn(M)$. Choose a maximal face $F$ with this property. Choose an
element $\overline{0}\neq \overline{a} \in \gp(M) \cap \R F/ \gp(M\cap
F)$ which is annihilated by $p$. We may assume that $a \in \relint
F$.

If follows from \ref{topological} that $H^{\dim F+1}_\mm(R)_{-a}
\neq 0$, since the poset $\Fc(M,a)$ consists all faces which contain
$F$ but not $F$ itself. But $H^{\dim F+1}_\mm(R)_{-pa} = 0$, because
$\Fc(pa)$ consists all faces which contain $F$ including $F$
itself. This poset is the face poset of a polytope and thus acyclic.
We have derived a contradiction to the assumption that $R$ is
$F$-injective, because the map
$H^{\dim F+1}_\mm(R)_{-a} \to H^{\dim F+1}_\mm(R)_{-pa}$ is not injective.
\end{proof}

\begin{cor}
Let $M \subseteq \Z^m$ be a positive affine monoid. Then
the following statements are equivalent:
\begin{enumerate}
\item
$M$ is normal;
\item
$R$ is $F$-split for every field $K$ of characteristic $p>0$;
\item
$R$ is $F$-pure for every field $K$ of characteristic $p>0$;
\item
$R$ is $F$-injective for every field $K$ of characteristic $p>0$.
\end{enumerate}
\end{cor}
\begin{proof}
It is easy to see that $M$ is normal if and only if $\gp(M) \cap \R F=
\gp( M \cap F)$ for all faces of $\cn(M)$. Thus (i) is equivalent to
(ii), (iii) and (iv) by \ref{charp1} and \ref{primeprop}.
\end{proof}

In a sense, it is inessential for the results of this section that
$K$ has characteristic $p$. In order to obtain variants that are
valid for every field $K$, one has to replace the Frobenius
endomorphism by the natural inclusion $K[pM]\to K[M]$.

%------------------------------------------------------------------------
%
%
%
%------------------------------------------------------------------------
\section{Examples and Counterexamples}\label{ExandCEx}
In this section we present various examples and counterexamples
related to the results of this paper. We choose a field $K$.

We saw in \ref{cn_simp} that $R=K[M]$ is Cohen--Macaulay for a
positive seminormal monoid $M$ of rank $d \leq 3$. Since a
Cohen--Macaulay ring always satisfies $(S_2)$ one could hope that
seminormality and $(S_2)$ already imply the Cohen--Macaulay property
of $R$. This is not the case as the following example shows.

\begin{ex}
\label{exsemi} For this example and the following one we fix some
notation. Let $P$ be a $3$-dimensional pyramid with a square base embedded into
$\R^4$ in degree 1. For example let $P$ be the convex hull of the
vertices given by
\begin{gather*}
m_0=(0,0,1,1),\ m_1=(-1,1,0,1),\ m_2=(-1,-1,0,1),\\
m_3=(1,-1,0,1),\ m_4=(1,1,0,1).
\end{gather*}
Figure \ref{pyramid} shows projections of the pyramid onto its base.
\begin{figure}[bht]
\begin{pspicture}(0,0)(2,2.3)
 \pspolygon(0,0)(2,0)(2,2)(0,2)
 \pspolygon[fillstyle=solid,fillcolor=light](0,0)(1,1)(0,2)
 \pspolygon[fillstyle=solid,fillcolor=light](2,0)(1,1)(2,2)
 \rput(0,0){\vertex}
 \rput(2,0){\vertex}
 \rput(2,2){\vertex}
 \rput(0,2){\vertex}
 \rput(1,1){\vertex}
 \rput(-0.2,-0.2){$m_2$}
 \rput(-0.2,2.2){$m_1$}
 \rput(2.2,-0.2){$m_3$}
 \rput(2.2,2.2){$m_4$}
 \rput(1.4,1.0){$m_0$}
\end{pspicture}
\qquad\qquad\qquad
\begin{pspicture}(0,0)(2,2)
 \pspolygon(0,0)(2,0)(2,2)(0,2)
 \pspolygon[fillstyle=solid,fillcolor=light](0,0)(1,1)(0,2)
 \pspolygon(2,0)(1,1)(2,2)
 \rput(2,0){\vertex}
 \rput(2,2){\vertex}
 \rput(0,2){\vertex}
 \rput(1,1){\vertex}
 \rput(-0.2,-0.2){$m_2$}
 \rput(-0.2,2.2){$m_1$}
 \rput(2.2,-0.2){$m_3$}
 \rput(2.2,2.2){$m_4$}
 \rput(1.4,1.0){$m_0$}
\end{pspicture}
\caption{}\label{pyramid}
\end{figure}
Let $C$ be the cone generated by $P$, so that $P$ is a cross section
of $C$. The facets of $C$ are the cones
\begin{gather*}
F_0=\cn (m_1,m_2,m_3,m_4),\ F_1=\cn (m_0,m_1,m_2),\ F_2=\cn(m_0,m_2,m_3),\\
F_3=\cn (m_0,m_3,m_4),\ F_4=\cn (m_0,m_1,m_4).
\end{gather*}
Let $M$ be the monoid generated by all integer points of {\em even
degree} in the facets $F_1$ and $F_3$ and all their faces, and all
remaining integer points in the interior of all other faces of $C$
including $C$.
(The facets $F_1$ and $F_3$ have been shaded in the left diagram in
Figure \ref{pyramid}.) Thus $M$ is positive, $\cn(M)=C$ and $M$ is
not normal. It follows from \ref{semicharacter} that $M$ is
seminormal and from \ref{s2cor} that $R$ satisfies $(S_2)$.

We claim that $R$ is not Cohen--Macaulay. Observe that all faces of
$C$ are simple except the face $\cn(m_0)$. Thus \ref{semihelper2}
and \ref{simplecrit} imply that $H^i_\mm(R)_{-a}$ can be nonzero
only for some $a \in \cn(m_0)\cap \gp(M)$. Choose an arbitrary $a$
in the relative interior of $\cn(m_0)\cap \gp(M)$ of odd degree. The
set $\Fc(M,a)$ introduced in \ref{topological} is
$$
\{  F_2, F_4, C\}
$$
and we see that the complex $\mathcal{C}^{\pnt}(M,a)$ has cohomology
in cohomological degree $3$. Hence $H^{3}_\mm(R)_{-a} \neq 0$ and
therefore $R$ is not Cohen--Macaulay.

The reader may check that $R/\overline\omega$ is Cohen--Macaulay, but
both $K[F_2\cap M]$ and $K[F_4\cap M]$ have nonzero third local
cohomology in degree $-a$.
\end{ex}

\begin{rem}\label{best_poss}
Example \ref{exsemi} can be generalized in the following way: if $C$
is not a simple cone, then there exists a seminormal affine monoid
$M$ with $C=\cn(M)$ such that $K[M]$ satisfies ($S_2$), but is not
Cohen--Macaulay for any field $K$.
\end{rem}

Next we consider the question whether the Cohen--Macaulay property
or the $(S_2)$ property are inherited by face projections. A
counterexample to this claim is already given in \cite[Example
2.2]{HORO74}. We can modify \ref{exsemi} a little bit to get the
same result for seminormal monoids.

\begin{ex}
\label{excmfaces} With the same notation as in \ref{exsemi} let $C$
be the cone over the pyramid $P$ with facets $F_0,\dots,F_4$. Now
let $M$ be the monoid generated by all integer points of {\em even
degree} in the facet $F_1$ and all its faces (as indicated in the
right diagram in Figure \ref{pyramid}), and all remaining integer
points in the interior of all other faces of $C$. Thus $M$ is
positive, $\cn(M)=C$ and $M$ is not normal. It still follows from
\ref{semicharacter} that $M$ is seminormal and by \ref{s2cor} that
$R$ satisfies $(S_2)$. Since all proper faces of $C$ except
$\cn(m_0)$ are simple, we only have to check the vanishing of the
local cohomology for points in the in $-\relint \cn(m_0)$. Let $a
\in \relint \cn(m_0)$. If $a$ has even degree, then
$$
\Fc(M,a)=\{  F \text{ face of } C: \cn(m_0) \subseteq F\}.
$$
If $a$ has odd degree, then
$$
\Fc(M,a)=\{F_2\cap F_3, F_3\cap F_4, F_2, F_3, F_4,  C \}.
$$
In any case, we can check that the complex
$\mathcal{C}^{\pnt}(M,a)$ is acyclic and therefore
$H^{i}_\mm(R)=0$ for $i<\rank M$. Thus $R$ is Cohen--Macaulay and
must satisfy $(S_2)$.

But $K[M\cap F_3]$ has only depth $1$, as can be seen from a similar
discussion as for $R$. So it does not satisfy $(S_2)$. Hence neither
the Cohen--Macaulay property, nor $(S_2)$ are inherited by face
projections of seminormal monoid rings.
\end{ex}

Let $\Delta$ be a simplicial complex contained in the simplex
$\Sigma$  with vertex set $V$. We consider the dual simplex
$\Sigma^*$ whose facets correspond bijectively to the vertices $v\in
V$ of $\Delta$. Next we erect the pyramid $\Pi$ over $\Sigma^*$ with
apex $t$. Then the faces of $\Pi$ that contain $t$ are in bijective
correspondence with the faces $G$ of $\Sigma$:
$$
F\in\Sigma \quad\longleftrightarrow\quad F^*\in\Sigma^*
\quad\longleftrightarrow\quad \tilde F=F^*\circ \{t\}
$$
where $\circ$ denotes the join. Observe that this correspondence
reverses the partial order by inclusion. Choose a realization of
$\Pi$ as a rational polytope, also denoted by $\Pi$.

Next we plane off those faces of $\Pi$ that correspond to the
minimal non-faces of $\Delta$. For such a non-face $G$ we choose a
support hyperplane $H$ of $\Pi$ with $\Pi\cap H=\tilde G$. Moving
this hyperplane by a sufficiently small rational displacement
towards the interior of $\Pi$, and intersecting $\Pi$ with the
positive halfspace of the displaced parallel hyperplane $H'$ we
obtain a polytope $\Pi'$ such that exactly the faces $F$ of $\Sigma$
that do not contain $G$ are preserved in $\Pi'$: $\tilde F\cap
\Pi'\neq\emptyset\iff F\not\supset G$.

Repeating this construction for each minimal non-face of $\Delta$ we
finally reach a polytope $\Pi''$ in which exactly the faces $\tilde
F$, $F\in\Delta$, have survived in the sense that $F'= \tilde F\cap
\Pi''$ is a non-empty face of $\Pi''$.

Moreover, the only
facets of $\Pi''$ containing $F'$ are the facets $\{v\}'$
corresponding to the vertices $v\in F$. On the other hand, every
face $E$ of $\Pi''$ that is not of the form $F'$ is contained in at
least one ``new'' facet of $\Pi''$ created by the planing of $\Pi$.

Note that $\Pi$ is a simplex and therefore a simple polytope. The
process by which we have created $\Pi''$ does not
destroy simplicity if the displacements of the hyperplanes are
sufficiently small and ``generic''.

\begin{figure}[hbt]
\footnotesize
\begin{pspicture}(0,-0.5)(1,0)
 \rput(0,0){\vertex}
 \rput(1,0){\vertex}
 \rput(0.5,-0.3){$\Delta$}
 \rput(0,0.3){$v_1$}
 \rput(1,0.3){$v_2$}
\end{pspicture}
\quad\quad
\begin{pspicture}(0,-0.5)(1,0)
 \rput(0,0){\vertex}
 \rput(1,0){\vertex}
 \psline(0,0)(1,0)
 \rput(0.5,-0.3){$\Sigma^*$}
\end{pspicture}
\qquad\quad
\begin{pspicture}(0,0)(2,2)
 \rput(0,0){\vertex}
 \rput(2,0){\vertex}
 \rput(1,1.5){\vertex}
 \pspolygon(0,0)(2,0)(1,1.5)
 \rput(1,0.3){$\Pi$}
 \rput(-0.0,0.8){$\widetilde{\{v_1\}}$}
 \rput(2.0,0.8){$\widetilde{\{v_2\}}$}
 \rput(1,1.8){$\{v_1,v_2\}^\sim$}
\end{pspicture}
\qquad\qquad
\begin{pspicture}(0,0)(2,1.5)
 \rput(0,0){\vertex}
 \rput(2,0){\vertex}
 \rput(0.667,1){\vertex}
 \rput(1.333,1){\vertex}
 \pspolygon(0,0)(2,0)(1.333,1)(0.667,1)
 \psline[linestyle=dotted](-0.3,1)(2.3,1)
 \psline[linestyle=dotted](1.333,1)(1,1.5)(0.667,1)
 \rput(1,0.3){$\Pi''$}
 \rput(-0.2,0.6){$\{v_1\}'$}
 \rput(2.2,0.6){$\{v_2\}'$}
\end{pspicture}
\qquad\quad
\begin{pspicture}(0,0)(3,1.5)
 \rput(0,0){\vertex}
 \rput(2,0){\vertex}
 \rput(0.667,1){\vertex}
 \rput(1.333,1){\vertex}
 \rput(3,0.5){\vertex}
 \pspolygon(0,0)(2,0)(1.333,1)(0.667,1)
 \psline(0,0)(3,0.5)
 \psline(2,0)(3,0.5)
 \psline(0.667,1)(3,0.5)
 \psline(1.333,1)(3,0.5)
 \rput(1,0.5){$\Gamma$}
 \end{pspicture}
\caption{Planing off a face and the construction of
$\Gamma$}\label{planing}
\end{figure}
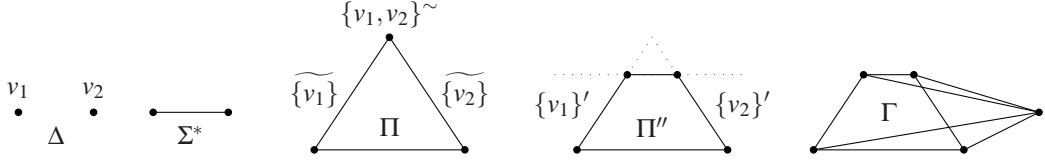

Set $d=\dim \Pi''+2$ and embed $\Pi''$ into
$\R^{d-2}\times\{0\}\subset\R^{d-1}$. Then let $\Gamma$ be the
pyramid over $\Pi''$ with apex $v=(0,\dots,0,1)$. The construction
of $\Gamma$ that leads to the pyramid of Example \ref{exsemi} is
illustrated in Figure \ref{planing}.

Note that all faces of $\Gamma$, except $\{v\}$, are simple.
($\{v\}$ is simple only if $\Delta=\Sigma$, or equivalently,
$\Pi''=\Pi$.)

In the last step we embed $\Gamma$ into $\R^{d-1}\times\{1\}\subset
\R^d$ by the assignment $x\mapsto (x,1)$ and choose the cone
$C=\R_+\Gamma$. It has dimension $d$. The point $v$ (in $\R^d$) has
the coordinates $(0,\dots,0,1,1)$. Therefore it has value $1$ under
the linear form $\deg:\R^d\to \R$, $\deg(y)=y_{d}$.

Set $L=\{a\in\Z^d:\deg(a)\equiv 0\pod 2\}$. To each facet $F$ of $C$
we assign the lattice
$$
L_F=\begin{cases} \R F\cap\Z^d,& F=\R_+\{v\}'\text{ for some
}v\in V,\\
\R F\cap L,&\text{else}.
\end{cases}
$$
Finally, we let $M$ be the monoid formed by all $a\in C\cap\Z^d$ such
that $a\in L_F$ for all facets $F$ of $C$ containing $a$. Clearly
$M$ is seminormal. Moreover, with the notation of Corollary
\ref{s2cor}, we have $M=M'$, since we have restricted the lattice
facet-wise, and thus $K[M]$ satisfies $(S_2)$ for all fields $K$.

Let $a\in\overline M=\gp(M)\cap C$. (By construction we have
$\gp(M)=\Z^d$.) If the face $F$ of $C$ with $a\in\relint(F)$ is
different from $\R_+v$, then it is a simple face of $C$, and
$H_\mm^i(K[M])_{-a}=0$ for all $i<n$ by Proposition
\ref{simplecrit}. If $F=\R_+v$ and $\deg(a)\equiv 0\pod 2$, then we
arrive at the same conclusion by Corollary \ref{van_cor}. However,
if $F=\R_+v$ and $\deg(a)\equiv 1\pod 2$, then the poset $\Fc(M,a)$ is
isomorphic to the dual of $\Delta$ (as a poset). Hence the cochain
complex $\mathcal{C}^\bullet(M,a)$ is isomorphic to the chain
complex of $\Delta$ (up to shift). We have
\begin{equation}\label{HDelta}
H_\mm^i(K[M])_{-a}=\tilde H_{d-i-1}(\Delta;K),\quad i=0,\dots,d.
\end{equation}

\begin{thm}
Let $\Delta$ be a simplicial complex and $K$ a field. Then there
exists a seminormal monoid $M$ of rank $d$ with $M=M'$ and such that
$R=K[M]$ has the following properties:
\begin{enumerate}
\item For every $a\in\overline M$
\begin{enumerate}
\item $H_\mm^i(R)_{-a}=0$ for all $i<d$, or
\item $H_\mm^\bullet(R)_{-a}$ is given by \eqref{HDelta}.
\end{enumerate}
\item Moreover, case \textup{(b)} holds true for at least one $a\in\overline M$.
\item The following are equivalent:
\begin{enumerate}
\item $\Delta$ is acyclic over $K$;
\item $R$ is Cohen-Macaulay.
\end{enumerate}
\end{enumerate}
\end{thm}

If we choose $\Delta$ as a triangulation of the real projective
plane, then we obtain a monoid algebra $K[M]$ which is
Cohen-Macaulay if and only if $\chara K\neq 2$.

\end{document}